\documentclass{amsart}
\address{M.S.R.I. \\
1000 Centennial Dr.\\
Berkeley, CA 94720 }
\email{jbryan@msri.org}

\title{Seiberg-Witten Theory and $\znums /2^{p}$ actions on spin 4-manifolds}
\author{Jim Bryan}
\date{}

\usepackage{amsmath,amsfonts}
\usepackage{amscd}

\newtheorem{thm}{Theorem}[section]
\newtheorem{cor}[thm]{Corollary}
\newtheorem{lem}[thm]{Lemma}
\newtheorem{prop}[thm]{Proposition}

\theoremstyle{definition}

\theoremstyle{remark}
\newtheorem{rem}{Remark}[section]
\newtheorem{example}{Example}[section]

\newcommand{\cnums} {{\mathbf C}}          
\newcommand{\rnums} {{\mathbf R}}		
\newcommand{\znums} {{\mathbf Z}}		
\newcommand{\qnums} {{\mathbf Q}}		
\newcommand{\cpbar} {\overline{\mathbf {CP}}^2}

\newcommand{\dirac}{\partial
\hspace{-5pt}\raisebox{1.5pt}{/}\hspace{-1pt}}

\newcommand{\Ker}{\operatorname{Ker}}
\newcommand{\End}{\operatorname{End}}
\newcommand{\Tr}{\operatorname{tr}}

\newcommand{\Coker}{\operatorname{Coker}}

\newcommand{\ie}{{\em i.e. }}

\newcommand{\til}[1]{\Tilde{#1}}

\newcommand{\1}{{{\mathchoice {\rm 1\mskip-4mu l} {\rm 1\mskip-4mu l}
{\rm 1\mskip-4.5mu l} {\rm 1\mskip-5mu l}}}}
\newcommand{\Zp}{\znums /2^{p}}
\newcommand{\Zi}{\znums /2^{i}}
\newcommand{\eleveneights}{``11/8-th's'' }
\newcommand{\teneights}{``10/8-th's'' }

\newcommand{\embed}{\hookrightarrow}
\newcommand{\CP}{\cnums \mathbf{P}^{2}}
\newcommand{\pin}{\operatorname{Pin}{(2)}}

\newcommand{\rootminus}{\sqrt{-1}}
\newcommand{\Ztwoq}{(\znums /2)^{q}}

\begin{document}

\begin{abstract}
Furuta's \teneights theorem gives a bound on the magnitude of the signature of a
smooth spin 4-manifold in terms of the second Betti number. We show that in
the presence of a $\Zp $ action, his bound can be strengthened. As
applications, we give new genus bounds on classes with divisibility and we
give a classification of involutions on  rational cohomology $K3$'s.

We utilize the action of $Pin(2)\til{\times }\Zp $ on the Seiberg-Witten
moduli space. Our techniques also provide a simplification of the proof of
Furuta's theorem.
\end{abstract}
\thanks{Supported by M.S.R.I. postdoctoral fellowship.}
\maketitle
\markboth{$\Zp $ actions on spin manifolds}
{$\Zp $ actions on spin manifolds}
\renewcommand{\sectionmark}[1]{}

\section{Introduction}\label{sec: intro}
In early 1995, Furuta \cite{Furuta} proved that if $X$ is a smooth,
compact, connected
spin 4-manifold with non-zero signature $\sigma (X)$, then

\begin{equation}\label{eq: 108ths}
\frac{5}{4}|\sigma (X)|+2 \leq b_{2}(X).
\end{equation}

This estimate has been dubbed the \teneights theorem in comparison with the
\eleveneights conjecture which predicts the following bound:
$$
\frac{11}{8}|\sigma (X)|\leq b_{2}(X).
$$
The inequality (\ref{eq: 108ths}) follows by a surgery argument from the
non-positive signature, $b_{1}(X)=0$ case:
\begin{thm}[Furuta]\label{thm: Furuta's}
Let $X$ be a smooth spin 4-manifold with $b_{1}(X)=0$ with non-positive
signature. Let $k=-\sigma (X)/16$ and $m=b_{2}^{+}(X)$. Then,
$$
2k+1\leq m
$$
if $m\neq 0$.
\end{thm}

The main result of this paper improves the above bound by $p$ under the
assumption that $X$ has a $\Zp $ action satisfying some non-degeneracy
conditions (the analogues of the condition $m\neq 0$ in the above theorem).

A $\Zp $ action is called a
{\em spin action} if the generator of the action $\tau :X\to X$ lifts to
the spin bundle $\hat{\tau }:P_{Spin}\to P_{Spin}$. Such an action is of
{\em even type} if $\hat{\tau }$ has order $2^{p}$ and is of {\em odd type}
if $\hat{\tau }$ has order $2^{p+1}$. Using \cite{At-Bo}, it easy to
determine if an action is spin and whether it is of even or odd type.

We assume throughout that $X$ is a smooth spin 4-manifold with $b_{1}(X)=0$
and oriented so that the signature is non-positive. We will continue with
the notation  $k=-\sigma (X)/16$ and $m=b_{2}^{+}(X)$. Our main results are
the following:
\begin{thm}\label{thm: thmA}
Suppose that $\tau:X\to X$ generates a smooth
$\Zp $ action that is spin and of odd type.  Let $X_{i}$ denote the
quotient of $X$ by $\Zi \subset \Zp$. Then
$$
2k+1+p\leq m
$$
if $m\neq 2k+b_{2}^{+}(X_{1})$ and $b_{2}^{+}(X_{i})\neq b_{2}^{+}(X_{j}) >0$
for $i\neq j$.
\end{thm}

For even actions, the easiest results to state are for involutions:
\begin{thm}\label{thm: thmB}
Suppose that $\sigma :X\to X$ is a smooth involution of even type. Then
$$
2k+2\leq m
$$
if  $m\neq b_{2}^{+}(X/\sigma )>0$.
\end{thm}

This is a special case of the following theorem:
\begin{thm}\label{thm: thmD}
Suppose that $\sigma _{1},\ldots,\sigma _{q}:X\to X$ are smooth involutions
of even type generating an action of $\Ztwoq$. Then
$$
2k+1+q\leq m
$$
if $m\neq b_{2}^{+}(X/g)$ for any non-trivial element $g\in \Ztwoq $ and
$b_{2}^{+}(X/\Ztwoq )\neq 0$.
\end{thm}

\begin{rem}\label{rem: non-deg conditions can also be in terms of fixed set}
The non-degeneracy conditions in Theorems \ref{thm: thmA}, \ref{thm:
thmB}, and \ref{thm: thmD} can be restated in terms of the fixed point set
although for our applications, the conditions on the quotient are more
convenient.
\end{rem}

We also include the following
\begin{thm}\label{thm: thmC}
Suppose that $\tau:X\to X$ generates a smooth
$\Zp $ action that is spin and of either type. Then $b_{2}^{+}(X/\tau )\neq
0$ if $k>0$.
\end{thm}

One can apply Theorem \ref{thm: thmA}  to a cover of a four
manifold branched along a smoothly embedded surface. The inequality can be
used to get a bound on the genus of the embedded surface (see Theorem
\ref{thm: genus bound in the general case}). One special case is the following:
\begin{thm}\label{thm: genus bound using double cover}
Let $M$ be a smooth, compact, oriented, simply connected 4-manifold (not
necessarily spin) with $b_{2}^{+}(M)>1$ and let $\Sigma \embed M$ be a
smooth embedding of a genus $g$ surface. Suppose that the homology class
defined by $\Sigma $ is divisible by 2 and that $[\Sigma] /2\equiv w_{2}(M)
\mod 2$. Then
$$
g\geq \frac{5}{4}\left(\frac{[\Sigma] ^{2}}{4} -\sigma (M)\right)-b_{2}(M)+2.
$$
\end{thm}

This bound is typically weaker than the adjunction inequality from the
Seiberg-Witten invariants (see \cite{M-S-T}), but the bound applies to
manifolds even with zero Seiberg-Witten invariants. For examples where the
bound is sharp and new we have (see also example \ref{example}):
\begin{cor}\label{cor: gmin for ((4,4),6) and (2,6)}
The minimal genus of an embedded surface in $\CP \# \CP $ representing the
class $(6,2)$ is 10. The minimal genus of an embedded surface in
$S^{2}\times S^{2} \# \CP $ representing the class $((4,4),6)$ is 19.
\end{cor}
\begin{proof}
In each case the lower bound of Theorem \ref{thm: genus bound using double
cover} is realized by the connected sum of smooth algebraic curves in each
factor.
\end{proof}

As another application of the theorems we give a ``classification'' of
involutions on rational cohomology $K3$'s.
\begin{thm}\label{thm: classification of involutions on Q-K3's}
Let $\sigma :X\to X$ be a spin involution of a rational cohomology $K3$
(\ie $b_{1}(X)=0$ and $Q_{X}\cong Q_{K3}$ ). If $\sigma $ is of even
type then it has exactly 8 isolated fixed points and $b_{2}^{+}(X/\sigma
)=3$; if $\sigma $ is of odd type then $b_{2}^{+}(X/\sigma )=1$.
\end{thm}

This theorem recovers as a special case a theorem of Donaldson concerning
involutions on the $K3$ (\cite{D-K} Cor. 9.1.4) and is related to  a
theorem of Ruberman  \cite{Ru}. We also remark that both possibilities in
the theorem actually occur.

The proof of Theorems \ref{thm: thmA}, \ref{thm: thmB},  \ref{thm: thmD},
and \ref{thm: thmC}
uses Furuta's technique of ``finite dimensional approximation'' for the
Seiberg-Witten
moduli spaces to reduce the problem to algebraic topology. The main
innovation of our technique is our approach to the equivariant $K$-theory.
In particular, we do not need the Adam's operations in equivariant
$K$-theory and can thus simplify that part of Furuta's proof.
In section  \ref{sec: S-W theory} we introduce the equations and use
Furuta's technique to study the moduli space; in section \ref{sec: equiv
K-theory} we use equivariant $K$-theory and representation theory to study
the $G$-equivariant properties of the moduli space and prove the main
theorems; section \ref{sec: applications} is devoted to applications,
primarily genus bounds obtained by branched covers and our classification
of involutions on rational cohomology $K3$'s.

We would like to thank R. Cohen, B. Gompf, C. Gordon, D. Kotschick, P.
Kronheimer, T. Lawson, N. C. Leung, G. Mat\'{\i}c, P. Oszvath, D. Ruberman,
R. Stern, P. Teichner, and R. Wentworth
for helpful conversations, and we would like to thank M.S.R.I. for
providing an excellent research environment where this work was carried
out.

\section{Seiberg-Witten Theory}\label{sec: S-W theory}
In this section we use Furuta's ``finite dimensional approximation''
technique to study the Seiberg-Witten moduli space in the presence of our
$\Zp $ symmetry. The goal of the section is to use the Seiberg-Witten
solutions to produce a certain $G$-equivariant map between spheres. This
map will then be studied by algebraic topology in section
\ref{sec: equiv K-theory} to produce a proof of the main theorems.

Let $X$ be a smooth, compact, connected, spin 4-manifold with $b_{1}(X)=0$
and fix an orientation so that $\sigma (X)\leq 0$. By Rochlin's theorem,
$\sigma (X) $ is divisible by 16 and so let $k=-\sigma (X)/16$ and let
$m=b_{2}^{+}(X)$. Let  $\tau :X\to  X$ be an orientation preserving
diffeomorphism generating a $\Zp $ action and fix an invariant Riemannian
metric $g$. We assume that $\tau $ is a spin action. By definition this 
means that $\tau^{*}(\wp )-\wp $ is zero as an element of $H^{1}(X;\znums
/2)$ where $\wp $ is the spin structure (the difference of two spin
structures is naturally an element of $H^{1}(X;\znums /2)$). There is then
a lift $\hat{\tau }$ of $\tau $ to the spin bundle:
$$
\begin{CD}
P_{\text{Spin}}&@>{\hat{\tau }}>>&P_{\text{Spin}}\\
@VVV&	&	@VVV\\
P_{\text{SO(4)}}&@>{d{\tau }}>>&P_{\text{SO(4)}}\\
@VVV&	&	@VVV\\
X&	@>{\tau }>>&	X.
\end{CD}
$$
The order of $\hat{\tau }$ is either $2^{p}$ or $2^{p+1}$ depending on
whether $(\hat{\tau })^{2^{p}}$ is the identity or the non-trivial deck
transformation of the double cover $P_{\text{Spin}}\to P_{\text{SO(4)}}$.
We say that $\tau $ is of even type if $\hat{\tau }$ has order $2^{p}$ and
of odd type if $\hat{\tau }$ has order $2^{p+1}$. A lemma of Atiyah and
Bott makes it easy to determine the type of $\tau $:
The $2^{p-1}$-th iteration of $\tau $ is an
involution $\sigma $ of $X$ and will have a fixed point set consisting of
manifolds of constant even dimension. If the fixed point set of $\sigma $
consists of 
points (or all of $X$), then $\tau $ is of even type; if the fixed point set
is dimension 2, then $\tau $ is of odd type. In the case where $\sigma $
acts freely, then $\tau $ is of even or odd type depending on whether the
quotient $X/\sigma $ is spin or not (c.f. \cite{At-Bo} or \cite{Br-Gok}).

On a spin manifold, the Seiberg-Witten equations for the trivial
spin$^{\cnums }$ structure take on a somewhat special form. In
\cite{Furuta}, Furuta gives a concise and elegant description of the
equations. We give a slightly different description to avoid repetition and
to make the exposition more closely match standard notation.

Let $S=S^{+}\oplus S^{-}$ denote the decomposition of the spinor bundle
into the positive and negative spinor bundles. Clifford multiplication
induces an isomorphism
$$\rho :\Lambda^{*}_{\cnums }\to \End _{\cnums }(S) $$
between the bundle of  complex valued forms and endomorphisms of $S$. Let
$\dirac :\Gamma (S^{+})\to \Gamma (S^{-})$ be the Dirac operator. The
Seiberg-Witten equations are for a pair $(a,\phi )\in \Omega
^{1}(X,\rootminus \rnums
) \times \Gamma (S^{+})$ and they are
\begin{eqnarray*}
\dirac \phi +\rho (a)\phi &	=&	0,\\
\rho (d^{+}a)-\phi \otimes \phi ^{*}+\tfrac{1}{2}|\phi |^{2}\1 &=&0,\\
d^{*}a&	=&	0.
\end{eqnarray*}

This system of equations is elliptic as written---the last equation defines
a slice for the (based) gauge group of the usual equations. The image of
$\rootminus \Lambda ^{2}_{+}$ under $\rho $ is the tracefree, hermitian
endomorphisms of
$S^{+}$ which we denote by $\rootminus \mathfrak{su}(S^{+})$. (We remark
that by
taking a slightly different gauge fixing condition, one could combine the last
two equations in the single equation $\rho (d^{+}a+d^{*}a)=\phi \otimes
\phi ^{*}.$)

Following Furuta's notation, we regard the equations as the zero set of a map
$$
D+Q:V\to W'
$$
where $V$ is the $L_{2}^{4}$-completion of $\Gamma (\rootminus \Lambda
^{1}\oplus
S^{+})$ and $W'$ is the $L^{3}_{2}$ completion of $\Gamma (S^{-}\oplus
\rootminus \mathfrak{su}(S^{+})\oplus \rootminus \Lambda ^{0})$ and $D$ and
$Q$ are the linear and quadratic pieces of the equations respectively, \ie 
\begin{eqnarray*}
D(a,\phi )&	=&	(\dirac \phi ,\rho (d^{+}a),d^{*}a),\\
Q(a,\phi )&	=&	(\rho (a)\phi ,\phi \otimes \phi
^{*}-\tfrac{1}{2}|\phi |^{2}\1 ,0).
\end{eqnarray*}

The image of $D+Q$ is $L^{2}$-orthogonal to the constant functions in
$\rootminus \Omega ^{0}\subset W'$. We define $W$ to be the orthogonal
complement of the constant functions in $W'$ and consider $D+Q$ as defined
on $W$:
$$
D+Q:V\to W.
$$

We now wish to determine the group of symmetries of the equations. As in
Furuta, we will see that $\pin$ acts and, as one would expect, there are
additional symmetries arising from $\tau $. The symmetry group turns out to
depend on whether $\tau $ is of even or odd type.

Define  $\pin \subset \text{SU}(2)$ to be the centralizer of $S^{1}\subset
\text{SU}(2)$. Regarding $\text{SU(2)}$ as the group of unit quaternions
and taking $S^{1}$ to be elements of the form
$e^{\rootminus \theta }$, $\pin $  then consists of elements of the form
$e^{\rootminus \theta }$ or $e^{\rootminus \theta }J$ (our quaternions are
spanned by $\langle 1,\rootminus ,J,\rootminus J \rangle$). We define the
action of $\pin $ on $V$
and $W$ as follows: Since $S^{+}$  and $S^{-}$ are $\text{SU(2)}$ bundles,
$\pin $ naturally acts on $\Gamma (S^{\pm })$ by multiplication on the left.
$\znums /2$ acts on
$\Gamma (\Lambda ^{*}_{\cnums })$ by multiplication by $\pm 1$ and this pulls
back to an action of $\pin $ by the natural map $\pin \to \znums /2$. A
calculation shows that this pullback also describes the induced action of
$\pin $ on $\rootminus \mathfrak{su}(S^{+})$. Both $D$ and $Q$ are then
seen to be
$\pin $ equivariant maps. Note that the action of $S^{1}\subset \pin $ is
the ordinary action of the constant gauge transformations.

The isometry $\tau $ acts on $V$ and $W$ by pull back by $d\tau $ and
$\hat{\tau }$; $D$ and $Q$ are equivariant with respect to this action.
If $\tau $ is of even type, then it induces an action of $\Zp $ on $V$ and
$W$; if $\tau $ is of odd type then it induces an action of $\znums
/2^{p+1}$. In the even case the symmetry group is thus
$$
G_{ev}\equiv \pin \times \znums /2^{p}.
$$
In the odd case $\pin \times \znums /2^{p+1}$ acts but we see that the
diagonal $\znums /2$ subgroup acts trivially: the $2^{p}$-th iteration of
$\hat{\tau } $ acts on the spinors by the action induced by the non-trivial deck
transformation; this is the same as the action of the constant gauge
transformation $-1\subset S^{1}\subset \pin $. Thus the symmetry group for
the odd case is
$$
G_{odd}\equiv \frac{\pin \times \znums /2^{p+1}}{\znums /2}.
$$
Throughout the sequel we write $G$ to mean $G_{odd}$ or $G_{ev}$.

To state Furuta's finite dimensional approximation theorem we define
$V_{\lambda }$ for any $\lambda \in \rnums $ to be the subspace of $V$
spanned by the eigenspaces of $D^{*}D$ with eigenvalues less than or equal
to $\lambda $. Similarly, define $W_{\lambda }$ using $DD^{*}$ and write
$p _{\lambda }$ for the $L^{2}$-orthogonal projection from $W$
 to $W_{\lambda }$. Define
$$
D_{\lambda }+Q_{\lambda }:V_{\lambda }\to W_{\lambda }
$$
by the restriction of $D+p_{\lambda }Q$ to $V_{\lambda }$.

Since the eigenspaces of $D^{*}D$ and $DD^{*}$ are $G$ invariant and
$p_{\lambda }$ is a $G$-map, $D_{\lambda }+Q_{\lambda }$ is a $G$
equivariant map between finite dimensional $G$ representations and it is an
approximation to $D+Q$ in the following sense:
\begin{lem}[Lemma 3.4 of \cite{Furuta}]\label{lem: D+Qlambda-1(0) is compact}
There exists an $R\in \rnums $ such that for all $\lambda >R$ the inverse
image of zero $(D_{\lambda }+Q_{\lambda })^{-1}(0)$ is compact.
\end{lem}

The lemma allows Furuta to use $D_{\lambda }+Q_{\lambda }$ to construct a
$G$-equivariant map on disks preserving boundaries
$$
f_{\lambda }:(BV_{\lambda ,\cnums },SV_{\lambda ,\cnums })\to (BW_{\lambda
,\cnums },SW_{\lambda ,\cnums })  .
$$
Here $V_{\lambda ,\cnums }=V_{\lambda }\otimes \cnums $, $BV_{\lambda
,\cnums }$ is homotopic to a ball in $V_{\lambda ,\cnums }$, and
$SV_{\lambda ,\cnums }$ is the boundary of $BV_{\lambda ,\cnums }$ with
similar definitions for $W_{\lambda ,\cnums }$, $BW_{\lambda ,\cnums }$,
and $SW_{\lambda ,\cnums }$.

The virtual $G$-representation $[V_{\lambda }]-[W_{\lambda }]\in R(G)$ is
the $G$-index of $D$ and can be determined by the $G$ index theorem and is
independent of $\lambda $. We discuss its computation in the next section.

\section{Equivariant $K$-theory}\label{sec: equiv K-theory}

In this section we use equivariant $K$-theory to deduce restrictions on the
map $f_{\lambda }:BV_{\lambda ,\cnums }\to BW_{\lambda ,\cnums }$.
Combining this with the index theorem determining $[V_{\lambda }]-[W_{\lambda
}]\in R(G)$ we will prove the main theorems. Our $K$-theoretic techniques
avoid Furuta's use of the equivariant Adam's operations and thus also
provide a simplification in the proof of his \teneights theorem.

\subsection{The representation ring of $G_{ev}$ and
$G_{odd}$.}\label{subsec: R(G)}
We write $R(\Gamma )$ for the complex representation ring of a compact Lie
group $\Gamma $ and we write direct sum and tensor product of
representations additively and
multiplicatively respectively. The representation ring $R(\Zp )$ of $\Zp $
is isomorphic to the group ring $\znums (\Zp )$ and is generated by the
standard one dimensional representation $\zeta$ . We write $1$
 for the trivial representation, so for example $\zeta ^{2^{p}}=1$, and as
a $\znums $ module, $R(\Zp )$ is generated by $1,\zeta ,\ldots,\zeta 
^{2^{p}-1}.$

The group $\pin $ has one non-trivial one dimensional
representation which we denote $\til{1}$ given by pulling back the
non-trivial $\znums /2$ representation by the map $\pin \to \znums /2.$ It
has a countable series of 2 dimensional irreducible representations
$h_{1},h_{2},\ldots$. The representation $h_{1}$, which we sometimes 
write as $h$, is the restriction of the standard representation of $SU(2)$
to $\pin \subset SU(2)$. The representations $h_{i}$ can be obtained using the
relation $h_{i}h_{j}=h_{i+j}+h_{|i-j|}$ where by convention $h_{0}$ denotes
$1+\til{1}$. Note that $\til{1}\cdot h_{i}=h_{i}$. Let $\theta $ denote the
standard 1 dimensional representation of $S^{1}$ so that $R(S^{1})\cong
\znums [\theta ,\theta ^{-1}]$. It is easy to see that $h_{i}$ restricts to
$\theta ^{i}+\theta ^{-i}$ as an $S^{1}$ representation.

Since $G_{ev}=\pin \times \Zp $, the representation ring is just the tensor
product $R(G_{ev})=R(\pin )\otimes R(\Zp )$. We can thus write a general
element $\beta \in R(G_{ev})$ as
$$
\beta =\beta _{0}(\zeta )1+\til{\beta }_{0}(\zeta )\til{1}+\sum_{i=1}^{\infty
}\beta _{i}(\zeta )h_{i}
$$
where $\beta _{0}$, $\til{\beta}_{0}$, $\beta _{1}$,\dots  are degree $2^{p}-1$
polynomials in $\zeta $ and all but a finite number of the $\beta _{i}$'s
are 0.

The irreducible representations of $G_{odd}=(\pin \times \znums
/2^{p+1})/(\znums /2)$ are the representations of $\pin \times \znums
/2^{p+1}$ that are invariant under $\znums /2$. To avoid confusion we will
use $\xi $ for the generator of $R(\znums /2^{p+1})$. The $\pin \times
\znums /2^{p+1}$ representations $\xi ^{i}$ and $h_{i}$ are non-trivial
restricted to $\znums /2$ if and only if $i$ is odd and the representation
$\til{1}$ is trivial restricted to $\znums /2$.
The subring $R(G_{odd})\subset R(\pin \times \znums /2^{p+1})$ is therefore
generated by $1,\til{1},\xi ^{2},$ and $\xi ^{i}h_{j}$ where $i\equiv j
\pmod 2$. We write a general element $\beta \in R(G_{odd})$ as
$$
\beta =\beta _{0}(\xi )+\til{\beta }_{0}(\xi )+\sum_{i=1}^{\infty }\beta
_{i}(\xi )h_{i}
$$
where now $\beta _{0},\til{\beta }_{0},$ and $\beta _{2i}$ are {\em even}
polynomials of degree $2^{p+1}-2$ and the $\beta _{2i+1}$'s are {\em odd}
polynomials of degree $2^{p+1}-1$. In summary we have:
\begin{thm}\label{thm: the relations of R(G)}
The ring $R(G_{ev})$ is generated by $1,$ $\til{1},$ $h_{1},$
$h_{2},\ldots,$ and
$\zeta $ with the relations $\zeta ^{2^{p}}=1$ and
$h_{i}h_{j}=h_{i+j}+h_{|i-j|}$, where $h_{0}=1+\til{1}$.

The ring $R(G_{odd})$ is generated by $1,$ $\til{1},$ $\xi ^{2},$ $h_{2},$
$h_{4},\ldots,$ and
$\xi h_{1},$ $\xi h_{3},\ldots $ with the relations $\xi  ^{2^{p+1}}=1$ and
$h_{i}h_{j}=h_{i+j}+h_{|i-j|}$, where $h_{0}=1+\til{1}$.

For $G=G_{odd}$ or $G_{ev}$ the restriction map $R(G)\to R(S^{1})\cong
\znums [\theta ,\theta ^{-1}]$ is given by $\til{1}\mapsto 1,$ $ \xi
\mapsto 1$ (or $\zeta\mapsto 1$), and $h_{i}\mapsto \theta ^{i}+\theta ^{-i}$.
\end{thm}

\subsection{The index of $D$.}\label{subsec: index of D}

The virtual representation $[V_{\lambda ,\cnums }]-[W_{\lambda ,\cnums
}]\in R(G)$ is the same as $\operatorname{Ind}(D)=[\Ker D]-[\Coker D]$.
Furuta determines $\operatorname{Ind}(D)$ as a $\pin $
representation; denoting the restriction map $r:R(G)\to R(\pin )$, Furuta shows
$$
r(\operatorname{Ind}(D))=2kh-m\til{1}
$$
where $k=-\sigma(X)/16$ and $m=b^{+}_{2}(X)$. Thus
$\operatorname{Ind}(D)=sh-t\til{1}$ where
$s$ and $t$ are polynomials in $\xi $ or $\zeta $ such that $s(1)=2k$ and
$t(1)=m$. In the case of $G=G_{ev}$ we write
$$
\begin{array}{ccc}
s(\zeta )=\sum_{i=1}^{2^{p}}s_{i}\zeta ^{i}&\text{and}&t(\zeta
)=\sum_{i=1}^{2^{p}}t_{i}\zeta ^{i}
\end{array}
$$
and for $G=G_{odd}$ we have
$$
\begin{array}{ccc}
s(\xi )=\sum_{i=1}^{2^{p}}s_{i}\xi ^{2i-1}&\text{and}&t(\xi
)=\sum_{i=1}^{2^{p}}t_{i}\xi ^{2i} .
\end{array}
$$
To compute the coefficients $t_{i}$ and $s_{i}$ one can use the $\Zp $-index
theorem. We will only need information about the $t_{i} $ coefficients.

From the definition of the operator $D$ we see that the polynomial
$\sum_{i=1}^{2^{p}}t_{i}\zeta ^{i}$ is the honest  $\Zp $
representation
$$\Coker ((d^{*},d^{+}):\Omega ^{1}\to \Omega
^{0}_{\bot}\oplus \Omega ^{2}_{+})$$
where $\Omega ^{0}_{\bot}$ denotes the $L^{2}$-orthogonal complement of the
constant functions (this is true for both the odd and even cases).
The constant coefficient, $t_{2^{p}}$ is the dimension of the invariant part
of $H^{2}_{+}(X)$ which is just $b_{2}^{+}(X/\tau )$. More generally,
$b_{2}^{+}$ of the quotients of $X$ by subgroups of $\Zp $ are given by the
various sums of the $t_{i}$'s. The dimension of the subspace of
$\sum_{i=1}^{2^{p}}t_{i}\xi ^{i}$ invariant under $\znums /2^{j}\subset
\znums /2^{p}$ is given by
$$
\sum_{i\equiv 0 \bmod 2^{j}}t_{i}=b_{2}^{+}(X_{j})
$$
where $X_{j}$ is the quotient of $X$ by $\znums /2^{j}$ and $X_{0}=X$ by
convention.

We summarize the discussion in the following
\begin{thm}\label{thm: index of D}
The index of $D$ is given by
$$
[V_{\lambda ,\cnums }]-[W_{\lambda ,\cnums }]=sh-t\til{1}\in R(G)
$$
where if $G=G_{ev}$ then
$$
\begin{array}{ccc}
s(\zeta )=\sum_{i=1}^{2^{p}}s_{i}\zeta ^{i}&\text{and}&t(\zeta
)=\sum_{i=1}^{2^{p}}t_{i}\zeta ^{i}
\end{array}
$$
and if $G=G_{odd}$ then
$$
\begin{array}{ccc}
s(\xi )=\sum_{i=1}^{2^{p}}s_{i}\xi ^{2i-1}&\text{and}&t(\xi
)=\sum_{i=1}^{2^{p}}t_{i}\xi ^{2i} .
\end{array}
$$

For either $G_{ev} $ or $G_{odd}$ we have
\begin{eqnarray*}
\sum_{i=1}^{2^{p}}s_{i}&	=&	2k\\
\sum_{i\equiv 0\bmod 2^{j}}t_{i}&=&b^{+}_{2}(X_{j}).
\end{eqnarray*}
\end{thm}

\subsection{The Thom isomorphism and a character formula for the
$K$-theoretic degree.}\label{subsec: Thom iso and char formula}
The Thom isomorphism theorem in equivariant $K$-theory for a general
compact Lie group is a deep theorem proved using elliptic operators
\cite{At-elliptic} . The subsequent character formula of this section
uses only
elementary properties of the Bott class. We follow tom Dieck \cite{tD} pgs.
254--255 for this discussion.

Let $V$ and $W$ be complex $\Gamma $ representations for some compact Lie
group $\Gamma $. Let $BV$ and $BW$ denote balls in $V$ and $W$ and let
$f:BV\to BW$ be a $\Gamma $-map preserving the boundaries $SV$ and $SW$.
$K_{\Gamma }(V)$ is by definition $K_{\Gamma }(BV,SV)$ and by the
equivariant Thom isomorphism theorem, $K_{\Gamma }(V)$ is a free $R(\Gamma
)$ module with generator the Bott class $\lambda (V)$. Applying the
$K$-theory functor to $f$ we get a map
$$
f^{*}:K_{\Gamma }(W)\to K_{\Gamma }(V)
$$
which defines a unique element $\alpha _{f}\in R(\Gamma )$ by the equation
$f^{*}(\lambda (W))=\alpha _{f}\cdot \lambda (V)$. The element $\alpha
_{f}$ is called the {\em $K$-theoretic degree} of $f$.

Let $V_{g}$ and $W_{g}$ denote the subspaces of $V$ and $W$ fixed by an
element $g\in \Gamma $ and let $V_{g}^{\bot}$ and $W_{g}^{\bot}$ be the
orthogonal complements. Let $f^{g}:V_{g}\to W_{g}$ be the restriction of
$f$ (well defined because of equivariance) and let $d(f^{g})$ denote the
ordinary topological degree of $f^{g}$ (by definition, $d(f^{g})=0$ if
$\dim V_{g}\neq \dim
W_{g}$). For any $\beta \in R(\Gamma )$ let $\lambda _{-1}\beta $ denote
the alternating sum $\sum_{}(-1)^{i}\lambda ^{i}\beta $ of exterior powers.

T. tom Dieck proves a character formula for the degree $\alpha _{f}$:
\begin{thm}\label{thm: char formula}
Let $f:BV\to BW$ be a $\Gamma $-map preserving boundaries and let $\alpha
_{f}\in R(\Gamma )$ be the $K$-theory degree. Then
$$
\Tr _{g}(\alpha _{f})=d(f^{g})\Tr _{g}(\lambda _{-1}(W_{g}^{\bot}-V_{g}^{\bot}))
$$
where $\Tr _{g}$ is the trace of the action of an element $g\in \Gamma $.
\end{thm}

This formula will be especially useful in the case where $\dim W_{g}\neq
\dim V_{g}$ so that $d(f^{g})=0$.

We also recall here that $\lambda
_{-1}(\sum_{i}a_{i}r_{i})=\prod_{i}(\lambda _{-1}r_{i})^{a_{i}}$ and that
for a one dimensional representation $r$, we have $\lambda _{-1}r=(1-r)$. A
two dimensional representation such as $h$ has $\lambda _{-1}h=(1-h+\Lambda
^{2} h)$. In this case, since $h$ comes from an $SU(2)$ representation,
$\Lambda ^{2}h=\det h=1$ so $\lambda _{-1}h=(2-h)$.

All the proofs in the following subsections proceed by using the character
formula to
examine the $K$-theory degree $\alpha _{f_{\lambda }}$ of the map
$f_{\lambda }:BV_{\lambda ,\cnums }\to BW_{\lambda ,\cnums }$ coming from
the Seiberg-Witten equations. We will abbreviate $\alpha _{f_{\lambda }}$ as
just $\alpha $ and $V_{\lambda ,\cnums }$ and $W_{\lambda, \cnums }$ as
just $V$ and $W$ and we will use the following elements of $G$. Let $\phi \in
S^{1}\subset \pin \subset G$ be an element generating a dense subgroup of
$S^{1}$; let $\eta \in \Zp $ and $\nu \in \znums /2^{p+1}$ be generators
and recall that there is the element $J\in \pin $ coming from the
quaternions. Note that the action of $J$ on $h$ has two invariant subspaces
on which $J$ acts by multiplication with $\rootminus $ and $-\rootminus $.

\subsection{Proof of Furuta's theorem (Theorem \ref{thm: Furuta's})}

Consider $\alpha =\alpha _{f_{\lambda }}\in R(\pin )$; it has the form
$$
\alpha =\alpha _{0}+\til{\alpha }_{0}\til{1}+\sum_{i=1}^{\infty }\alpha
_{i}h_{i} .
$$
Since $\phi $ acts non-trivially on $h$ and trivially on $\til{1}$, we have
that $\dim V_{\phi }\neq \dim W_{\phi }$ as long as $m>0$. The character
formula then gives
$$
\Tr _{\phi }(\alpha )=0=\alpha _{0}+\til{\alpha }_{0}+\sum_{i=1}^{\infty
}\alpha _{i}(\phi ^{i}+\phi ^{-i})
$$
so $\alpha _{0}=-\til{\alpha }_{0}$ and $\alpha _{i}=0$ for $i\geq 1$.

Since $J$ acts non-trivially on both $h$ and $\til{1}$, $\dim V_{J}=\dim
W_{J}=0$ so $d(f^{J})=1$ and the character formula says
\begin{eqnarray*}
\Tr _{J}(\alpha )&	=&	\Tr _{J}(\lambda _{-1}(m\til{1}-2kh))\\
&	=&	\Tr _{J}((1-\til{1})^{m}(2-h)^{-2k})\\
&	=&	2^{m-2k}
\end{eqnarray*}
using $\Tr _{J}h=0$ and $\Tr _{J}\til{1}=-1$. On the other hand $\Tr
_{J}(\alpha )=\Tr _{J}(\alpha _{0}(1-\til{1}))=2\alpha _{0}$ so the degree is
$$
\alpha =2^{m-2k-1}(1-\til{1})
$$
and we can conclude $2k+1\leq m$. \qed

The proofs of Theorems \ref{thm: thmA}, \ref{thm: thmB}, and \ref{thm:
thmD} are a generalization
of the above proof.

\subsection{Proof of Theorem \ref{thm: thmB}}
In this case the group is
$G=G_{ev}=\pin \times \znums /2$ and
$$
V-W=(s_{1}\zeta +s_{2})h-(t_{1}\zeta +t_{2})\til{1}
$$
where $t_{1}+t_{2}=m$ and $s_{1}+s_{2}=2k$. The hypothesis $m\neq
b_{2}^{+}(X/\sigma )>0$ translates to $t_{1}+t_{2}\neq t_{2}>0$ which is
equivalent to both $t_{1}$ and $t_{2}$ being non-zero.

Both $\phi $ and $\phi \nu $ act non-trivially on $h$ and trivially on
$\til{1}$ so (since $t_{2}\neq 0$) we have $d(f^{\phi })=d(f^{\phi \nu
})=0$ so that $\Tr _{\phi }(\alpha )=\Tr _{\phi \nu }(\alpha )=0$.

The general form of $\alpha $ is
$$
\alpha =\alpha _{0}(\zeta )+\til{\alpha }_{0}(\zeta )\til{1}
+\sum_{i=1}^{\infty }\alpha _{i}(\zeta )h_{i}
$$
so the conditions $\Tr _{\phi }(\alpha )=\Tr _{\phi \nu }(\alpha )=0$ imply
that
\begin{eqnarray*}
\alpha _{0}(\pm 1)+\til{\alpha }_{0}(\pm 1)&	=&	0\\
\alpha _{i}(\pm 1)&	=&	0
\end{eqnarray*}
since $\alpha _{0},$ $\til{\alpha }_{0}$, and $\alpha _{i}$ are degree 1 in
$\zeta $ we see that $\alpha _{0}=-\til{\alpha }_{0}$ and $\alpha _{i}=0$ so
$$
\alpha =(\alpha _{0}^{1}\zeta +\alpha _{0}^{2})(1-\til{1}).
$$
Now $J\nu $ acts trivially on $\zeta \til{1}$ and non-trivially on $\zeta h$
and $h$ so since $t_{1}\neq 0$, $d(f^{J\nu })=0$. Thus we get
$$
\Tr _{J\nu }(\alpha )=0=(-\alpha _{0}^{1}+\alpha _{0}^{2})\cdot 2
$$
so that $\alpha _{0}^{1}=\alpha _{0}^{2}$.

Finally, as before, $J$ acts non-trivially on $W$ and $V$ so that $d(f^{J})=1$
and $\Tr _{J}(\alpha )=2^{m-2k}$ so it must be the case that
$$
\alpha = 2^{m-2k-2}(1+\zeta )(1-\til{1})
$$
and thus $2k+2\leq m$.\qed

\subsection{Proof of Theorem \ref{thm: thmD}}
We generalize the preceding proof. Now $G=\pin \times \Ztwoq $ and we
write $\zeta _{i}$ for the non-trivial representation of the $i$-th copy of
$\znums /2$. The index $\operatorname{Ind}(D)$ is given by
$$
V-W=s(\zeta _{1},\ldots,\zeta _{q})h-t(\zeta _{1},\ldots,\zeta _{q})\til{1}
$$
where $s$ and $t$ are polynomial functions in the variables $\zeta
_{1},\ldots,\zeta _{q}$ of multi-degree $(1,\ldots,1)$.

The hypothesis $m\neq b_{2}^{+}(X/g)$ then implies that the representation
$t(\zeta _{1},\ldots,\zeta _{q})$ contains a summand on which $g$ acts as
$-1$. Since $J$ acts on
$\til{1} $ by $-1$ we see that the representation $t(\zeta
_{1},\ldots,\zeta _{q})\til{1}$ has a positive dimensional subspace fixed
by $Jg$ for every non-trivial $g\in \Ztwoq $. On the other hand, $Jg$
always acts non-trivially on $s(\zeta _{1},\ldots,\zeta _{q})h$ so the
character formula gives us
$$
\Tr _{Jg}(\alpha )=0.
$$

Since $b_{2}^{+}(X/\Ztwoq )\neq 0$, the coefficient of the trivial
representation in $t(\zeta _{1},\ldots,\zeta _{q})$ is non-zero and so
$\phi g$ always fixes a non-trivial subspace of $t\til{1}$. On the
other hand, $\phi g$ always acts non-trivially on $sh$ so again the
character formula shows that
$$
\Tr _{\phi g}(\alpha )=0.
$$

The general form of $\alpha $ is
$$
\alpha =\alpha _{0}+\til{\alpha }_{0}\til{1}+\sum_{i=1}^{\infty }\alpha
_{i}h_{i}
$$
where $\alpha _{0}$, $\til{\alpha }_{0}$, and $\alpha _{i}$ are polynomial
functions in $\zeta _{1},\ldots,\zeta _{q}$ of multi-degree $(1,\ldots,1)$.

Now for any polynomial function $\beta (x_{1},\ldots,x_{q})$ of
multi-degree $(1,\ldots,1)$, if we know that $\beta
((-1)^{n_{1}},\ldots,(-1)^{n_{q}})=0$ for any arrangement of the signs,
then $\beta \equiv 0$. Thus the formula $\Tr _{\phi g}(\alpha )=0$ for all
$g\in \Ztwoq $ implies that $\alpha _{0}+\til{\alpha }_{0}\equiv 0$ and
$\alpha _{i}\equiv 0$. We can therefore write
$$
\alpha =\alpha _{0}(\zeta _{1},\ldots,\zeta _{q})(1-\til{1}).
$$

Since $J$ acts non-trivially on both $U$ and $V$ we can compute as before: 
$$\Tr _{J}(\alpha )=\Tr _{J}(\lambda _{-1}(W-V))=2^{m-2k}.$$
This equation, along with the $2^{q}-1$ equations $\Tr _{Jg}(\alpha )=0$
for $g\neq 1\in \Ztwoq $ give $2^{q} $ independent conditions on $\alpha
_{0}$ which determine it uniquely. It must be the following:
$$
\alpha =2^{m-2k-1-q}(1+\zeta _{1})(1+\zeta _{2})\cdots (1+\zeta
_{q})(1-\til{1})
$$
and thus $2k+1+q\leq m$.\qed

\subsection{Proof of Theorem \ref{thm: thmA}}
Recall that
$$
[V]-[W]=s(\xi )h-t(\xi )\til{1}\in R(G_{odd})
$$
with $s(\xi )=\sum_{i=1}^{2^{p}}s_{i}\xi ^{2i-1}$ and $t(\xi
)=\sum_{i=1}^{2^{p}}t_{i}\xi ^{2i}$. The $K$-theory degree $\alpha
=\alpha _{f_{\lambda }}$ of $f_{\lambda }$ has the form
$$
\alpha =\alpha _{0}(\xi )+\til{\alpha }_{0}(\xi )+\sum_{i=1}^{\infty
}\alpha _{i}(\xi )h_{i} .
$$
We compute the following characters of $\alpha $:
\begin{lem}\label{lem: chars of alpha}
$\Tr _{\phi \eta ^{j}}(\alpha )=0$ for $j=1,\ldots,2^{p}, \Tr _{J\eta
^{j}}(\alpha )=0$ for $j=1,\ldots,2^{p}-1$, and $\Tr _{J}(\alpha
)=2^{m-2k}$.
\end{lem}
{\em Proof:} We introduce the notation $(V)_{g}=\dim V_{g}$. By the character
formula, $\Tr _{\phi \eta ^{j}}(\alpha )=0$ if
$$
(V)_{\phi \eta^{j}}-(W)_{\phi \eta ^{j}}=(s(\xi )h-t(\xi )\til{1})_{\phi
\eta ^{j}}\neq 0.
$$
Since $\phi \eta ^{j}$ acts non-trivially on every $\xi ^{i}h$ and $\phi
\eta ^{j} $ acts trivially on $\til{1}$, we have
$$
(V)_{\phi \eta ^{j}}-(W)_{\phi \eta ^{j}}\leq -t_{2^{p}}.
$$
By the hypothesis in Theorem \ref{thm: thmA}, $t_{2^{p}}=b_{2}^{+}(X/\tau
)>0$ so $\Tr _{\phi \eta ^{j}}(\alpha )=0$.

To show $\Tr _{J\eta ^{j}}(\alpha )=0$ for $j=1,\ldots,2^{p}-1$ we need to
show $(V)_{J\eta ^{j}}-(W)_{J\eta ^{j}}\neq 0. $ The 2 dimensional
representation $h$ decomposes into two complex lines on which $J$ acts as
$\rootminus $ and $-\rootminus $, so
$$
(\xi ^{2i-1}h)_{J\eta ^{j}}=\begin{cases}
1&	\text{if $\eta ^{j} $ acts on $\xi ^{2i-1}$ by $\pm \rootminus $}\\
0&	\text{otherwise.}
\end{cases}
$$
$\eta ^{j}$ acts as $\pm \rootminus $ on $\xi ^{2i-1}$ if and only if $j$
and $i$ satisfy
$$
(2i-1)j\equiv \pm 2^{p-1}\bmod 2^{p+1}.
$$
If $j=2^{p-1} $ then the condition is satisfied for every $i$; if $j\neq
2^{p-1}$ the condition is never satisfied (divide both sides by the highest
power of 2 that divides $j$ to get an odd number on one side and an even on
the other). Thus
$$
(s(\xi )h)_{J\eta ^{j}}=\begin{cases}
2k&	\text{if $j=2^{p-1}$,}\\
0&	\text{if $j\neq 2^{p-1}$}.
\end{cases}
$$

Now $J\eta ^{j}$ acts on $\til{1}$ by $-1$ and so

$$
(\xi ^{2i}\til{1})_{J\eta ^{j}}=\begin{cases}
1&	\text{if $\eta ^{j}$ acts on $\xi ^{2i}$ by $-1$ and}\\
0&	\text{otherwise.}
\end{cases}
$$
Thus we see that
\begin{eqnarray*}
(t(\xi )\til{1})_{J\eta ^{j}}&	=&\sum_{2ij\equiv 2^{p}\bmod 2^{p+1}} t_{i}\\
&	=&	\sum_{ij\equiv 2^{p-1}\bmod 2^{p}}t_i.
\end{eqnarray*}

We want to see that the non-degeneracy conditions of the theorem imply
$(V-W)_{J\eta ^{j}}\neq 0$. The conditions are
equivalent to
$$
0<b_{2}^{+}(X_{p})<b_{2}^{+}(X_{p-1})<\cdots <b_{2}^{+}(X_{1})
$$
and
$$
b_{2}^{+}(X_{1})\neq m-2k
$$
which in terms of the $t_{i}$'s are
\begin{eqnarray*}
t_{2^{p}}&\neq &0\\
t_{2^{p-1}}&	\neq &	0\\
t_{2^{p-2}}+t_{3\cdot 2^{p-2}}&	\neq &	0\\
t_{2^{p-3}}+t_{3\cdot 2^{p-3}}+t_{5\cdot 2^{p-3}} +t_{7\cdot
2^{p-3}}&	\neq &	0\\
&	\vdots&	\\
t_{2}+t_{3\cdot 2}+t_{5\cdot 2}+\cdots &\neq &	0
\end{eqnarray*}
and
$$
\sum_{i\text{ even}}t_{i} \neq m-2k.
$$
Since $m=\sum t_{i}$ the last condition is the same as
$$
t_{1}+t_{3}+\cdots +t_{2^{p}-1}\neq 2k.
$$

From the previous discussion we have in the case $j\neq 2^{p-1}$
\begin{eqnarray*}
(V-W)_{J\eta ^{j}}&	=&	(s(\xi )h-t(\xi )\til{1})_{J\eta ^{j}} \\
&	=&	-\sum_{ij\equiv 2^{p-1}\bmod 2^{p}}t_{i}\\
&	=&	-(t_{2^{a}}+t_{3\cdot 2^{a}}+\cdots )\\
&	\neq &	0
\end{eqnarray*}
where $2^{p-1-a}$ is the largest power of 2 dividing $j$.

For the case $j=2^{p-1}$
\begin{eqnarray*}
(V-W)_{J\eta ^{j}}&	=&2k-\sum_{i2^{p-1}\equiv 2^{p-1}\bmod 2^{p}}t_{i} \\
&	=&	2k-\sum_{i\equiv 1\bmod 2}t_{i}\\
&	\neq &	0.
\end{eqnarray*}

To complete the lemma we compute $\Tr _{J}(\alpha )$. Once again, since $J$
acts non-trivially on $W$ and $V$,
$$
\Tr _{J}(\alpha )=\Tr _{J}\Lambda _{-1}(W-V)=2^{m-2k}.
$$

\begin{lem}\label{lem: polynomial lemma}
Let $\beta (\xi )$ be of the form $\sum_{i=1}^{2^{p}}\beta _{i}\xi ^{2i}$
or $\sum_{i=1}^{2^{p}}\beta _{i}\xi ^{2i-1}$ in $R(\znums /2^{p+1})$. If
$\Tr _{\eta ^{j}}\beta =0$ for $j=1,\ldots,2^{p}$ then $\beta \equiv 0$.
\end{lem}
{\em Proof:} Let $\underline{\beta }(x)=\sum_{i=1}^{2^{p}}\beta _{i}x^{i}$.
$\underline{\beta }$ is then a degree $2^{p}$ polynomial with roots at all
of the $2^{p}$-th roots of unity and at 0, so $\underline{\beta }$ is
identically 0.

The general form of $\alpha $ is
$$
\alpha =\alpha _{0}(\xi )+\til{\alpha }_{0}(\xi )\til{1}
+\sum_{i=1}^{\infty }\alpha _{i}(\xi )h_{i}
$$
so using the computation of Lemma \ref{lem: chars of alpha}
\begin{eqnarray*}
0&	=&	\Tr _{\phi \eta ^{j}}(\alpha )\\
&	=&\Tr _{\eta ^{j}}(\alpha _{0}+\til{\alpha }_{0})
+\sum_{i=1}^{\infty }\Tr _{\eta ^{j}}(\alpha _{i})(\phi ^{i}+\phi ^{-i})
\end{eqnarray*}
for $j=1,\ldots,2^{p}-1$.

Since each $\phi ^{i}$ term must vanish separately, Lemma \ref{lem:
polynomial lemma} immediately implies that $\alpha _{0}+\til{\alpha
}_{0}\equiv 0$ and $\alpha _{i}\equiv 0$. Thus $\alpha $ has the form $\alpha
_{0}(\xi )(1-\til{1})$ where $\alpha _{0}(\xi )=\sum_{i=1}^{2^{p}}\alpha
_{0}^{i} \xi ^{2i}$. Since $\Tr_{J\eta ^{j}}(\alpha )=\Tr_{\eta
^{j}}(\alpha _{0})\cdot 2=0$ for $j=1,\ldots,2^{p}-1$, the degree $2^{p}$
polynomial
$$
\underline{\alpha }_{0}(x)=\sum_{i=1}^{2^{p}}\alpha _{0}^{i}x^{i}
$$
has $2^{p}$ known roots, namely the $2^{p}-1$ non-trivial $2^{p}$-th roots
of unity and 0. Thus $\underline{\alpha
}_{0}(x)=\text{const.}\sum_{i=1}^{2^{p}} x^{i}$ and we can use $\Tr
_{J}(\alpha )=2^{m-2k}$ to determine the constant. Thus we can conclude
that
$$
\alpha =2^{m-2k-p-1}(1+\xi ^{2}+\xi ^{4}+\cdots +\xi
^{2^{p+1}-2})(1-\til{1}),
$$
and so $2k+1+p\leq m$ and the theorem is proved. \qed
\begin{rem}\label{rem: arbitrary 2-groups}
For a spin action of an arbitrary group $\Gamma $ of order $2^{p}$ on $X$
our methods should give the bound $2k+1+p\leq m$ as long as the action is
subject to some non-degeneracy conditions. The proof would proceed as all
the others, using the non-degeneracy conditions and the character formula
to guarantee $\Tr _{J g}(\alpha )=0$ for non-trivial $g\in \Gamma $ and
$\Tr _{\phi g}(\alpha )=0$ for all $g\in \Gamma $. Then one can show that
$$
\alpha =2^{m-2k-1-p}\cdot \rho \cdot (1-\til{1})
$$
where $\rho $ is the regular representation of $\Gamma $. Complications do
occur: in order to guarantee that $\Tr _{Jg}(\alpha )=0$ one needs to
incorporate information about the virtual representation $s$, \ie the
$\Gamma $-index of the Dirac operator. This can be computed via the index
theorem, but in general the non-degeneracy conditions would then involve
the fixed point set. A further complication occurs because the lift of the
action to the spin bundle can be complicated; some group elements may lift
with twice their original order while others may preserve their order.
Nevertheless, our techniques could be applied in a case by case basis if
the above issues are understood.
\end{rem}

\subsection{Proof of Theorem \ref{thm: thmC}}
This has a slightly different flavor than the previous proofs. The proof is
essentially the same for $G_{odd}$ and $G_{ev}$; we will use the $G_{ev}$
notation. The index has the form
$$
\operatorname{Ind}(D)=s(\zeta )h-t(\zeta )\til{1}
$$
and the hypothesis $k>0$ means that $s(\zeta )=\sum_{i=1}^{2^{p}}s_{i}\zeta
^{i}$ has at least one positive coefficient. From Theorem \ref{thm:
Furuta's}, $k>0$ also implies $m>0$. We will prove that $b_{2}^{+}(X/\tau
)\neq 0$ by contradiction. Suppose $b_{2}^{+}(X/\tau )=0$ so that the
constant coefficient of $t(\zeta )$ is 0, \ie $t_{2^{p}}=0$. Since $m>0$ we
have $(W-V)_{\phi }\neq 0$ so $\Tr _{\phi }(\alpha )=0$ and since
$t_{2^{p}}=0$, we know that $\phi \nu $ acts non-trivially on all of $V$
and $W$ so
\begin{eqnarray*}
\Tr _{\phi \nu }(\alpha )&=&\lambda _{-1}(W-V)\\
&=&\frac{(1-\nu )^{t_{1}}(1-\nu ^{2})^{t_{2}}\cdots (1-\nu
^{2^{p}-1})^{t_{2^{p}-1}}}{(1-\nu \phi )^{s_{1}}(1-\nu \phi
^{-1})^{s_{1}}\cdots }.
\end{eqnarray*}

Since at least one of the $s_{i}$'s is positive, the above expression has
arbitrarily high powers of $\phi $ in it. On the other hand
\begin{eqnarray*}
\Tr _{\phi \nu }(\alpha )&=&\Tr _{\phi \nu }\left(\alpha _{0}(\zeta
)+\til{\alpha }_{0}(\zeta )\til{1}+\sum_{i=1}^{\infty }\alpha _{i}(\zeta
)h_{i} \right)\\
&	=&\alpha _{0}(\nu )+\til{\alpha }_{0}(\nu )+\sum_{i=1}^{\infty
}\alpha _{i}(\nu )(\phi ^{i}+\phi ^{-i})
\end{eqnarray*}
has only finitely many non-zero $\alpha _{i}$ terms, which is our
contradiction. \qed

\section{Applications}\label{sec: applications}
\subsection{Genus bounds.}\label{subsec: genus bounds}
Our original motivation for this work is our application to genus bounds.
We refer the reader to our paper \cite{Br-Gok} for details of the set up.

Let $M$ be a smooth 4-manifold (not necessarily spin) and let $\Sigma
\embed M$ be a smoothly embedded surface representing a homology class
$a\in H_{2}(M;\znums )$. A basic question in 4-manifold topology asks for
the minimal genus $g_{min}(a)$ of a smoothly embedded surface representing
a given class $a$. In order to determine $g_{min}(a)$ one needs good lower
bounds on the genus of embedded surfaces and constructions realizing those
bounds.

The basic ideas we use go back to Hsiang-Szczarba
\cite{Hs-Sz}, Rochlin \cite{Ro}, and Kotschick-Mat\'{\i}c \cite{Ko-Ma} (see
also T. Lawson \cite{Lawson}).

Suppose the class $a$ is divisible by some number $d$, then one can
construct the $d$-fold branched cover $X\to M$ branched along $\Sigma $.
Under favorable hypotheses, $X$ will be spin and the covering
transformation is a spin action. The signature and Euler characteristic of
$X$ can be computed in terms of $\Sigma \cdot \Sigma $ and $g(\Sigma )$ and
so the bounds of Theorems \ref{thm: Furuta's} and \ref{thm: thmA} give
genus bounds.

The most straightforward implementation of this idea is for double
branched covers and results in Theorem \ref{thm: genus bound using double
cover}. We will prove the general result coming from $2^{p}$ covers in this
section.

\begin{prop}\label{prop: branched cover is spin}
Let $\Sigma \embed M$ be a smoothly embedded surface of genus $g$ in a
smooth, simply-connected, compact, oriented 4-manifold. Suppose that
$[\Sigma ]\in H_{2}(M)$ has the property that $2^{p}|[\Sigma ]$ and
$[\Sigma ]/2^{p}\equiv w_{2}(M)\bmod 2$.

Then there is a spin 4-manifold $X$ with a spin $\znums /2^{p}$ action
$\tau :X\to X$  of odd
type  such that $X/\tau =M$. Furthermore, $k=-\sigma (X)/16$ and
$m=b_{2}^{+}(X)$ are given by:
\begin{eqnarray*}
k&=&\frac{1}{16}\left(-2^{p}\sigma (M)+\frac{4^{p}-1}{3\cdot 2^{p}}[\Sigma
]\cdot [\Sigma ] \right),\\
m&=&2^{p}b_{2}^{+}(M)+(2^{p}-1)g-\frac{4^{p}-1}{6\cdot 2^{p}}[\Sigma ]\cdot
[\Sigma ].
\end{eqnarray*}
\end{prop}

{\em Proof:} This is in \cite{Br-Gok}; a sketch of the proof is the
following. The condition $2^{p}|[\Sigma ]$ allows one to construct the
$2^{p}$-fold branched cover $X\to M$ branched along $\Sigma $. One can
compute $w_{2}(X)$ using a formula of Brand \cite{Brand} and the condition
of $[\Sigma ]/2^{p}\equiv w_{2}(M)\bmod 2$ guarantees $w_{2}(X)$ vanishes. 
$\pi _{1}(X)$ is finite (see \cite{Ro}, \cite{Ko-Ma}) so $b_{1}(X)=b_{3}(X)=0$
and so $m$ and $k$ can be determined by $e(X)$ and $\sigma (X)$ which can
be computed by the $G$-index theorem. The covering action $\tau :X\to X$ is
automatically spin because $H_{1}(X;\znums /2)=0$ (see \cite{Ko-Ma}) and the
action is of odd type since the fixed point set is two dimensional. \qed

Here is our general genus bound (Theorem \ref{thm: genus bound using double
cover} is a special case):

\begin{thm}\label{thm: genus bound in the general case}
Let $\Sigma \embed M$ be a smoothly embedded surface of genus $g$ in a
smooth, simply-connected, compact, oriented 4-manifold. Suppose that
$[\Sigma ]\in H_{2}(M)$ has the property that $2^{p}|[\Sigma ]$ and
$[\Sigma ]/2^{p}\equiv w_{2}(M)\bmod 2$.

Suppose $b_{2}^{+}(M)>1$ and $g\neq [\Sigma
]^{2}(1+2^{2i-2p+1})/6-b_{2}^{+}(M)$ for
$i=1,\ldots,p-1 $. Then
\begin{equation}\label{eqn: genus bound}
g\geq \frac{1}{2^{p}-1}\left[\frac{5}{4}\left(\frac{4^{p}-1}{6\cdot
2^{p}}[\Sigma ]^{2}-2^{p-1}\sigma (M) \right)+1+p-2^{p-1}b_{2}(M) \right].
\end{equation}
\end{thm}
\begin{rem}
The conditions on $b_{2}^{+}(M)$ and $g$ are so that the non-degeneracy
conditions of Theorem \ref{thm: thmA} are met. The condition on $g$ is not
present for $p=1$. Although it is
aesthetically unpleasing to have restrictions on $g$ in the hypothesis, in
practice the conditions are easy to dispense with (see example
\ref{example} ).
Even without the conditions one still gets an
inequality coming from Furuta's theorem which is the same as the bound
(\ref{eqn: genus bound}) but without the $p$ term in the brackets. Our theorem
thus improves the Furuta bound  by $p/(2^p-1)$. Probably the best
applications occur in the case $p=1$ (as in Theorem \ref{thm: genus bound
using double cover} where we improve the Furuta bound by 1). However, we
will give examples where Theorem \ref{thm: genus bound in the general
case} provides sharp bounds unobtainable from the Furuta bound.
\end{rem}

\begin{example}\label{example}
Let $M=\# _{N}{\bf CP}^{2}$ with $N>1$ and consider $\Sigma \embed M$
representing the class $(4,\ldots,4)$. We can apply Theorem \ref{thm: genus
bound in the general case} with $p=2$. It tells us
$$
g\geq \tfrac{1}{3}(8N+3)
$$
if $g\neq 3N$. When $N\geq 6$ the condition on $g$ is no condition and for
$N=2,\ldots,5$ the theorem gives $g\geq 3N$ since if $g=3N-1$ it would
contradict $g\geq\tfrac{1}{3}(8N+3)$. Since the connected sum of the
algebraic representative in each factor has genus $3N$, the theorem gives a
sharp bound for $N=2,\ldots,5$.
\end{example}

{\em Proof of Theorem \ref{thm: genus bound in the general case} :}
Let $X\to M$ be the $2^{p}$-th branched cover and, to match the notation of
Theorem \ref{thm: thmA}, let $X_{i}$ be the $2^{p-i}$ branched cover. Then
$X_{0}=X$, $X_{p}=M$ and $X_{i}$ is the quotient of $X$ by $\znums
/2^{i}\subset \Zp $. Let $m_{i}=b^{+}_{2}(X_{i})$ which is computable by
essentially the same formula as in Proposition \ref{prop: branched cover is
spin}  
$$
m_{i}=2^{p-i}m_{p}+(2^{p-i}-1)g -\frac{4^{p-i}-1}{6\cdot 2^{p-i}}[\Sigma
]^{2}.
$$
Note that $m_{0}\geq m_{1}\geq\cdots \geq m_{p}>1$. The inequality of the
theorem is equivalent to $2k+1+p\leq m_{0}$ so the only way the theorem can
fail is if $2k+\delta =m_{0}$ for $\delta \in \{1,\ldots,p \}$ {\em and}
$m_{0}-2k=m_{1}$ or $m_{i}=m_{i-1}$ for $i=2,\ldots,p$. The hypothesis
$g\neq [\Sigma
]^{2}(1+2^{2i-2p+1})/6-b_{2}^{+}(M)$
is equivalent to $m_{i}\neq
m_{i+1}$ so we know that $m_{1}>m_{2}>\cdots
>m_{p}=b_{2}^{+}(M)$. This implies $m_{1}\geq p-1+b_{2}^{+}(M)$ which, 
combined with the hypothesis $b_{2}^{+}(M)>1$ gives $m_{1}>p$. Then we see
that we cannot have  $2k+\delta =m_{0}$ and $m_{0}-2k=m_{1}$  since that would
imply $\delta >p$.\qed

\subsection{Involutions on rational cohomology $K3$'s.}\label{subsec:
involutions on K3}

The proof of Theorem \ref{thm: classification of involutions on Q-K3's}
follows readily from our main theorems (Theorems \ref{thm: thmA}, \ref{thm:
thmB}, and \ref{thm: thmC}). Suppose $X$ is a rational cohomology $K3$ so
that $H^{*}(X;\qnums)\cong H^{*}(K3;\qnums)$ and suppose $\tau :X\to X$ is
a spin involution.

First suppose that $\tau $ is of odd type. Since the inequality of Theorem
\ref{thm: thmA} is violated, one of the non-degeneracy conditions must also
be violated and so $b_{2}^{+}(X/\tau )$ is either 1 or 0. Theorem \ref{thm:
thmC} shows that $b_{2}^{+}(X/\tau )$ is not 0 so it must be 1. This can
easily occur, for example $K3$ is a double branched cover over ${\bf
CP}^{2}$, $S^{2}\times S^{2}$, and $E(1)\cong {\bf CP}^{2}\# 9\cpbar $ and
it has a free involution covering the Enriques surface.

Now suppose that $\tau $ is even so it must violate the non-degeneracy
conditions of Theorem \ref{thm: thmB}. In conjunction with Theorem
\ref{thm: thmC} we see that $b_{2}^{+}(X/\tau )$ must be 3. Since the fixed
point set is at most points, the $G$-signature theorem tells us that
$$
\sigma (X/\tau )=\tfrac{1}{2}\sigma (X)=-8.
$$
If $N$ is then the number of fixed points, then the Lefschetz formula tells
us the Euler characteristic
$$
\chi (X/\tau )=\tfrac{1}{2}(\chi (X) +N)=12+N/2.
$$
Since $b_{1}(X/\tau )=0$ and we know $b_{2}^{+}(X/\tau )=3$, we can solve
the above equations to get $b_{2}^{-}(X/\tau )=11$ and $N=8$.

It was pointed out to us by P. Kronheimer that this does occur. The
construction is as follows. Let $\til{Y}$ be the $K3$ surface. There are 8
disjoint $-2$ spheres $S_{1},\ldots,S_{8}$ such that
$S=\sum_{i=1}^{8}S_{i}$ is divisible by 2.  Let $\til{X}$ be the double
branched cover of $\til{Y}$ branched along $S$. The preimage of $S$ in
$\til{X}$ is then 8 disjoint $(-1)$ spheres which we can blow down to
obtain $X$, a smooth manifold with an involution fixing 8 points covering
$Y$, the orbifold obtained by collapsing the $S_{i}$'s in $\til{Y}$:

$$
\begin{CD}
\til{X}@>>>X\\
@VVV@VVV\\
\til{Y}@>>>Y.
\end{CD}
$$
We will show that $X$ is a rational cohomology $K3$ (with a little more
work, one can see that $X\cong K3$) and $b_{2}^{+}(Y)=3$ so this is the
example we seek. From the signature theorem and Lefschetz formula we
compute that $\chi (\til{X})=2\chi (\til{Y})-\chi (S)=32$ and $\sigma
(\til{X})=2\sigma (\til{Y})-S^{2}=24$ and so, after blowing down, $\chi
(X)=24$  and $\sigma (X)=-16$. Since the preimage of $S$ is characteristic
in $\til{X}$ we know that $X$ is spin and since
$b_{1}(\til{X})=b_{1}(X)=0$, $X$ is a rational cohomology $K3$.

\begin{rem}
It is easy to construct a {\em homeomorphism} $\tau :X\to X $ generating an
involution with $b_{2}^{+}(X/\tau )=2$ (for example) but our classification
implies this will not be
smoothable, \ie  there is no smooth structure on $X$ so that this $\tau $
is a diffeomorphism.
\end{rem}


\end{document}